\documentclass[11pt,a4paper]{article}
\usepackage{amsmath}
\usepackage{amssymb}
\usepackage{amsthm}
\usepackage{fullpage}
\usepackage{graphicx}
\usepackage{hyperref}
\usepackage[capitalise,noabbrev]{cleveref}
\usepackage{booktabs}
\usepackage{multirow}
\usepackage{color}
\usepackage[sc]{mathpazo}
\usepackage[T1]{fontenc}
\usepackage[margin=10pt,font=small]{caption}

\newtheorem{thm}{Theorem}[section]
\newtheorem{prop}[thm]{Proposition}
\newtheorem{lem}[thm]{Lemma}
\newtheorem{coro}[thm]{Corollary}
\theoremstyle{remark}
\newtheorem{rmk}[thm]{Remark}
\newtheorem{conj}[thm]{Conjecture}

\newcommand{\tdef}[1]{\textcolor{blue}{\emph{#1}}}
\newcommand{\horiz}{\begin{center}\rule{0.3\textwidth}{0.5pt}\end{center}}

\newcommand{\defeq}{\mathrel{\mathop:}=}
\newcommand{\occ}{\operatorname{occ}}
\newcommand{\maxocc}{\operatorname{maxocc}}
\newcommand{\swentropy}{S_{\mathrm{sw}}}
\newcommand{\minswentropy}{\min\!\swentropy}
\newcommand{\expect}{\mathbb{E}}
\newcommand{\prob}{\mathbb{P}}
\newcommand{\naturals}{\mathbb{N}}
\newcommand{\logb}{\log_2}

\author{Wenjie Fang \\ Univ Gustave Eiffel, CNRS, LIGM, F-77454 Marne-la-Vallée, France}

\title{Maximal number of subword occurrences in a word}
\begin{document}

\maketitle

\abstract{We consider the number of occurrences of subwords (non-consecutive sub-sequences) in a given word. We first define the notion of subword entropy of a given word that measures the maximal number of occurrences among all possible subwords. We then give upper and lower bounds of minimal subword entropy for words of fixed length in a fixed alphabet, and also showing that minimal subword entropy per letter has a limit value. A better upper bound of minimal subword entropy for a binary alphabet is then given by looking at certain families of periodic words. We also give some conjectures based on experimental observations.}

\horiz

\section{Introduction} \label{sec:intro}

Enumeration problems concerning patterns have been rich sources of interesting combinatorics. The most famous examples are classes of permutations avoiding a given pattern. We refer readers to \cite{kitaev-pattern-survey,perm-pattern-survey} for an exposition of such results. In this article, we will consider enumeration about patterns in a word, which is in general easier than that for permutations.

There are two different widely used notions of patterns for words. The first notion is that of a \emph{factor}. A word $v$ occurs in another word $w$ \emph{as a factor} if there is a consecutive segment of $w$ equal to $v$. The second notion is that of a \emph{subword}. A word $v$ occurs in another word $w$ \emph{as a subword} if we can obtain $v$ by deleting letters in $w$. A factor of $w$ is always a subword of $w$, but not \emph{vice versa}. There are also other notions of patterns, such as the one in \cite{burstein-mansour} that generalizes both factors and subwords, but we will not discuss them here.

Unlike for permutations, the enumeration of classes of words avoiding a (set of) given subwords or factors is already known in the sense that, for a given subword or factor, we can express their avoidance in regular expressions, leading automatically (no pun intended) to the generating function of such classes, which is always rational and can be effectively computed \cite[Section~V.5]{flajolet}. There is also some work on counting words with a fixed number of occurrences of a given pattern, for example \cite{burstein-mansour,subseq-freq}. For the other end of the spectrum, the problem of maximal density of certain patterns in words is considered by Burstein, Hästö and Mansour in \cite{packing-pattern}. Readers are referred to the survey-book of Kitaev \cite{kitaev-pattern-survey} for more of such results. In general, such results are non-trivial, due to the possible overlap of patterns.

We may also consider all patterns that occur in a given word. For the notion of factors, this idea leads to the notion of factor complexity of a word $w$, first defined by Morse and Hedlund in \cite{morse-hedlund} and also called ``subword complexity'', which is a function $f_w$ such that $f_w(k)$ is the number of distinct factors in $w$ of length $k$. In \cite{subword-complexity-enum}, Gheorghiciuc and Ward studied the factor complexity of random words. We may also want to consider the number of occurrences of a given pattern in a word. The work of Flajolet, Szpankowski and Vallée \cite{hidden-word-stat} establishes a Gaussian limit law and large deviations for the number of occurrences of a given subword in a long random word, again by analyzing overlap of occurrences of subwords. The number of occurrences of a given pattern is of particular interest in algorithmics with applications in data mining, in which researchers propose algorithms finding patterns with large number of occurrences \cite{extract-subsequence} and study complexity of such problems \cite{complexity-frequent}.

In this article, we take a further step on enumeration problems on patterns by considering the number of occurrences of all subwords in a word. We are motivated by the following decision problem: given a word $w$ and a number $n$, is there another word $u$ that occurs at least $n$ times as a subword in $w$? To our knowledge, it is unknown whether this natural problem in combinatorics of words can be solved in polynomial time. As an attempt to study this decision problem, it is thus a reasonable first step to explore the structure of subword occurrences in a given word.

For a given word $w$ that permits frequent occurrences of some subword $w'$, we may see it as having some ``large space'' for such a subword, and we would like to measure the ``extend'' of such space, or from the opposite direction, the ``disorder'' generated by the possible different occurrences. To this end, we define a notion of \emph{subword entropy}, which measures the maximal number of times that any subword can occur in a given word. We delay its precise definition to later sections. We then look at the minimal subword entropy of all words of a given length $n$ in an alphabet of $k$ letters, denoted by $\minswentropy^{(k)}(n)$, as it is easy to find the ones with maximal subword entropy. Using the super-additivity of minimal subword entropy, we show that $\minswentropy^{(k)}(n) / n$ has a finite limit $L_k$. We then concentrate on the binary case, showing some upper bounds of $L_2$ by looking at certain families of periodic words, inspired by experimental data. As a by-product, we also show that, given two words $w$ and $v$, the generating function of the number of occurrences of $v^r$ in $w^m$ is rational.

The rest of this article is organized as follows. We first give necessary definitions in \Cref{sec:prelim}, then some basic results on subword occurrences and minimal subword entropy in \Cref{sec:basic}, including the proof of the existence of the limit $L_k$ of $\minswentropy^{(k)}(n) / n$, and bounds of $L_k$. Then in \Cref{sec:binary}, we focus on the case of binary alphabet, and shows a better upper bound of $L_2$ than the one given in \Cref{sec:basic}. We end in \Cref{sec:open} with a discussion on open problems partially inspired by experimental results obtained for the binary case.

\paragraph{Acknowledgments} We thank Stéphane Vialette for bringing the question of subword occurrences of a given word to our attention, and for giving an idea for \Cref{prop:maxocc-lowerbound}. We thank Jim Fill for pointing out an omission at $n=3$ in the experimental results, and also for inspiring discussions. We also thank the anonymous reviewer for the constructive comments, especially for pointing out an error in the original proof of \Cref{prop:subword-000111}.

\section{Preliminaries} \label{sec:prelim}

A \tdef{word} $w$ of length $n$ is a sequence $w = (w_1, \ldots, w_n)$ of elements in a finite set $\mathcal{A}$ called the \tdef{alphabet}. We denote by $|w|$ the length of $w$, and $|w|_a$ the number of letters $a$ in $w$. For two words $v, w$, their concatenation is denoted by $v \cdot w$. We also denote by $\epsilon$ the empty word of length $0$. A \tdef{run} in a word $w$ is a maximal consecutive segment in $w$ formed by only one letter in $\mathcal{A}$.
Given a word $w$, if there is another word $w' = (w'_1, \ldots, w'_k)$ such that there is some set $P = \{ p_1 < \cdots < p_k \}$ of integers from $1$ to $n$ satisfying $w_{p_j} = w'_j$ for all $1 \leq j \leq k$, then we say that $w'$ is a \tdef{subword} of $w$, and we call the set $P$ an \tdef{occurrence} of $w'$ in $w$. We denote by $\occ(w, w')$ the number of occurrences of $w'$ in $w$. For instance, for $w = 011001$ and $w' = 01$, there are $5$ occurrences of $w'$ in $w$, which are $\{1, 2\}, \{1, 3\}, \{1, 6\}, \{4, 6\}, \{5, 6\}$. When $w'$ is not a subword of $w$, we have $\occ(w, w') = 0$, and when $w' = \epsilon$, we have $\occ(w, \epsilon) = 1$. In the literature (see, \textit{e.g.}, \cite{binomial-subword}), the number of occurrences of $w'$ in $w$ is also regarded as a generalization of the binomial coefficients, and is thus denoted by $\binom{w}{w'}$. We do not use this notation here to avoid confusion with binomial coefficients.

It is easy to find words who have a subword with a large number of occurrences. For instance, with $w = a^n$ for some letter $a \in \mathcal{A}$, the subword $w' = a^{\lfloor n/2 \rfloor}$ appears $\binom{n}{\lfloor n/2 \rfloor} \sim \left(\frac{2}{\pi n}\right)^{1/2} 2^n$ times. It is more difficult to find words in which no subword occurs frequently. To quantify such intuition, we define the \tdef{maximal subword occurrences} $\maxocc(w)$ of a word $w$ to be the maximal value of $\occ(w, w')$, and subwords $w'$ reaching this value are called \tdef{most frequent subwords} of $w$. We note that a word $w$ may have several most frequent subwords. We then define the \tdef{subword entropy} $\swentropy(w)$ of $w$ in an alphabet of size $k$ by
\[
  \swentropy(w) \defeq \logb \maxocc(w).
\]
We note that this definition does not depend on the size of the alphabet, as subword occurrences are fundamentally about subsets of positions, and the size of the alphabet is implicit in the word $w$. Now, finding words in which no subword occurs frequently is to find words minimizing their subword entropy. We define the minimal subword entropy for words of length $n$ in an alphabet of size $k$ by
\[
  \minswentropy^{(k)}(n) \defeq \min_{w \in \mathcal{A}^n, |\mathcal{A}| = k} \swentropy(w).
\]

\section{Some basic results} \label{sec:basic}

We start with some simple properties of $\occ(w, u)$.

\begin{lem} \label{lem:basics-3}
  For words $w, w', u, u'$, we have $\occ(w \cdot w', u \cdot u') \geq \occ(w, u) \occ(w', u')$.
\end{lem}
\begin{proof}
  Let $P$ (resp. $P'$) be an occurrence of $u$ (resp. $u'$) in $w$ (resp. $w'$). The set $Q = P \cup \{ p' + |w| \mid p' \in P'\}$ is an occurrence of $u \cdot u'$ in $w \cdot w'$, and the map $(P, P') \mapsto Q$ is clearly injective.
\end{proof}

\begin{lem} \label{lem:basics-4}
  For a word $w$, it has a most frequent subword $u$ with $w_1 = u_1$ and $w_{|w|} = u_{|u|}$.
\end{lem}
\begin{proof}
  Let $v$ be a most frequent subword of $w$. If $v_1 \neq w_1$, then for all occurrences $P$ of $v$ in $w$, we have $1 \notin P$. Then, $\{1\} \cup P$ is an occurrence of $v' = w_1 \cdot v$, which is thus also a most frequent subword. Otherwise, we take $v' = v$. We repeat the same reasoning on $v'$ for the last letter of $w$ to obtain $u$.
\end{proof}

We now give simple upper and lower bounds for $\maxocc(w)$ for any word $w$.

\begin{prop} \label{prop:maxocc-upperbound}
  Given an alphabet $\mathcal{A}$ of size $k$ and $n \geq 1$, for any word $w \in \mathcal{A}^n$, we have $\maxocc(w) \leq \binom{n}{\lceil n/2 \rceil}$, and it is realized exactly by $w = a^n$ for any letter $a \in \mathcal{A}$.
\end{prop}
\begin{proof}
  For $w'$ of length $k$, as occurrences of $w'$ in $w$ are subsets of $\{1, \ldots, n\}$, we have $\occ(w, w') \leq \binom{n}{k} \leq \binom{n}{\lceil n/2 \rceil}$. It is clear that only words composed by the same letter reach this bound.
\end{proof}

\begin{prop} \label{prop:maxocc-lowerbound}
  Given an alphabet $\mathcal{A}$ of size $k$ and $n \geq 1$, for any word $w \in \mathcal{A}^n$, we have
  \[
    \ln \maxocc(w) \geq \max_{0 \leq \ell \leq n} \ln \left( \binom{n}{\ell} k^{-\ell} \right).
  \]
  The right-hand side is asymptotically maximized at $\ell = \lfloor \frac{n}{k+1} \rfloor$, giving the asymptotic value $n \ln(1+k^{-1}) - \frac{1}{2} (\ln n) + O(1)$.
\end{prop}
\begin{proof}
  Let $u$ be a uniformly chosen word of length $\ell$. We have
  \[
    \expect[\occ(w, u)] = \sum_{P \subseteq \{1,\ldots,n\}, |P|=\ell} \prob[u \text{ occurs in } w \text{ at positions } P] = \binom{n}{\ell} k^{-\ell}.
  \]
  The first equality is from linearity of expectation, and the second from the fact that $u$ is uniformly chosen at random, and the probability does not depend on $P$. Hence, there is some $u^*$ with $\occ(w, u^*) \geq \expect[\occ(w, u)]$, implying the non-asymptotic part of our claim.

  For the asymptotic part, take $\alpha = \ell / n$. Using Stirling's approximation, we have
  \[
    \ln \left( \binom{n}{\ell} k^{-\ell} \right) = n \bigg[ -\alpha \ln\alpha - (1-\alpha)\ln(1-\alpha) - \alpha\ln k \bigg] - \frac{1}{2} \ln n + O(1).
  \]
  The coefficient of $n$ above is maximized for $\alpha = (k+1)^{-1}$, with value $\ln(1+k^{-1})$. We thus have our claim on the asymptotic growth.
\end{proof}

\begin{coro} \label{coro:sw-bounds}
  There are constants $c_1, c_2$ such that, for all $n \in \naturals$ and $w \in \mathcal{A}^n$ with $|\mathcal{A}| = k \geq 2$, we have
  \[
    \logb(1 + k^{-1}) n - \frac{1}{2} \logb n + c_1 \leq \minswentropy^{(k)}(n) \leq \swentropy(w) \leq n - \frac{1}{2} \logb n + c_2.
  \]
\end{coro}
\begin{proof}
  The bounds on $\swentropy(w)$ result from combining \Cref{prop:maxocc-upperbound} and \Cref{prop:maxocc-lowerbound} with $\ln \binom{n}{\lceil n / 2 \rceil} = n \ln 2 - \frac{1}{2} \ln n + O(1)$. The bounds for $\minswentropy^{(k)}(n)$ then follows.
\end{proof}

We now show that there is a limit for $\minswentropy^{(k)}(n)/n$. To this end, we need the well-known Fekete's lemma \cite{fekete} for super-additive sequences.

\begin{lem}[\cite{fekete}] \label{lem:fekete}
  Suppose that a sequence $(g_n)_{n \geq 1}$ satisfies that, for all $n, m \geq 1$, we have $g_{n + m} \geq g_n + g_m$. Then, for $n \to +\infty$, the value of $g_n/n$ either tends to $+\infty$ or converges to some limit $L$.
\end{lem}

We first show that the function $\minswentropy^{(k)}(n)$ is super-additive.

\begin{prop} \label{prop:super-additivity}
  Given $k \geq 2$, for any $n, m \geq 1$, we have
  \[
    \minswentropy^{(k)}(n + m) \geq \minswentropy^{(k)}(n) + \minswentropy^{(k)}(m).
  \]
\end{prop}
\begin{proof}
  Let $w$ be a word of length $n + m$ achieving minimal subword entropy $\minswentropy^{(k)}(n + m)$. We write $w = w' \cdot w''$, with $|w'| = n$ and $|w''| = m$. Let $v'$ (resp. $v''$) be a most frequent subword of $w'$ (resp. $w''$). We have
  \[
    \maxocc(w) \geq \occ(w' \cdot w'', v' \cdot v'') \geq \occ(w', v') \occ(w'', v'') = \maxocc(w') \maxocc(w'').
  \]
  The first inequality is from the definition of $\maxocc$, the second from \Cref{lem:basics-3}, and the equality comes from the definition of $v'$ and $v''$. By the definition of $w$, we have
  \[
    \minswentropy^{(k)}(n + m) \geq \logb \maxocc(w') + \logb \maxocc(w'') \geq \minswentropy^{(k)}(n) + \minswentropy^{(k)}(m).
  \]
  The second inequality is from the definition of $\minswentropy^{(k)}$.
\end{proof}

\begin{thm} \label{thm:swentropy-limit}
  For any $k \geq 2$, the sequence $(\minswentropy^{(k)}(n) / n)_{n \geq 1}$ converges to a certain limit $L_k < +\infty$.
\end{thm}
\begin{proof}
  \Cref{prop:super-additivity} shows that $\minswentropy^{(k)}(n)$ is super-additive. We then apply \Cref{lem:fekete}, and as $\minswentropy^{(k)}(n)/n$ is bounded above by some constant according to \Cref{coro:sw-bounds}, we have the existence of the limit $L_k$ which is finite.
\end{proof}

With the existence of the limit $L_k$, we can use known values of $\minswentropy^{(k)}(n)$ to give lower bounds for $L_k$.

\begin{prop} \label{prop:sw-limit-lower-bound}
  Given $k \geq 2$, we have $L_k \geq \minswentropy^{(k)}(n) / n$ for all $n$.
\end{prop}
\begin{proof}
  By iterating \Cref{prop:super-additivity}, we have $\minswentropy^{(k)}(rn) \geq r \minswentropy^{(k)}(n)$ for all $r \geq 1$. Diving both sides by $rn$, it means that the limit $L_k$ of $\minswentropy^{(k)}(rn) / rn$ is also larger than $\minswentropy^{(k)}(n)/n$.
\end{proof}

From \Cref{coro:sw-bounds}, we know that
\[
  \logb(1 + k^{-1}) \leq L_k \leq 1.
\]
When $k \to \infty$, the lower bound is asymptotically $(\ln 2)^{-1} k^{-1}$, which tends to $0$, while the upper bound stays constant. The next natural step is to try to give better bounds for $L_k$, and eventually compute the precise value of $L_k$. However, it seems to be a formidable task.

\section{Better upper bound for binary alphabet} \label{sec:binary}

After the general basic results given in \Cref{sec:basic}, we will focus hereinafter on the case of binary alphabet $\mathcal{A} = \{0, 1\}$. In this case, the bounds in \Cref{coro:sw-bounds} become $\logb(3/2) \leq L_2 \leq 1$ for the limit $L_2$ in \Cref{thm:swentropy-limit}. The gap between the two bounds are significant, as $\logb(3/2) \approx 0.585$. We now consider three families of periodic words with a small value of maximal subword occurrences, inspired and explored by computer experimentation (see \Cref{sec:open}). They are $(01)^m$, $(0011)^m$ and $(000111)^m$, all giving better upper bounds of $L_2$.

In the following, we present our results in the order of $(0011)^m$, $(01)^m$ and $(000111)^m$, as the first family gives the best upper bound, and the last one is the most complicated to analyze.

\subsection{Subword occurrences of $(0011)^m$}

We start by a structural result on most frequent subword of a word of the form $(0011)^m$.

\begin{prop} \label{prop:subword-0011}
  For $w = (0011)^m$, there is a most frequent subword $w'$ of the form $(01)^r$.
\end{prop}
\begin{proof}
  Take a most frequent subword $u$ of $w$ of length $\ell$. Suppose that $u$ has the form $u = s \cdot 00 \cdot t$. We take $u^{(1)} = s \cdot 010 \cdot t$ and $u^{(2)} = s \cdot 0 \cdot t$. Let $P = \{p_1, \ldots, p_\ell\}$ be an occurrence of $u$ in $w$, and we suppose that the $00$ occurs at $p_i, p_{i+1}$. Let $\mathcal{P}$ be the set of occurrences of $u$ in $w$, which is divided into $\mathcal{P} = \mathcal{P}_1 \cup \mathcal{P}_2$, where $\mathcal{P}_1$ contains those $P$'s with $p_i + 1 \neq p_{i+1}$, while $\mathcal{P}_2$ contains those with $p_i + 1 = p_{i+1}$. For any $P \in \mathcal{P}_1$, the two $0$'s occur in different runs, meaning that there is at least one run $11$ in between. This leads to at least two choices for the extra $1$ added in $u^{(1)}$. Therefore, $\occ(w, u^{(1)}) \geq 2 |\mathcal{P}_1|$. For any $P \in \mathcal{P}_2$, the two $0$'s occur in the same run, meaning that replacing them by a single $0$ leaves us two choices. We thus have $\occ(w, u^{(2)}) \geq 2 |\mathcal{P}_2|$, meaning that $\occ(w, u^{(1)}) + \occ(w, u^{(2)}) \geq 2 |\mathcal{P}| = 2 \maxocc(w)$. We deduce that at least one of the $u^{(j)}$'s satisfies $\occ(w, u^{(j)}) = \maxocc(w)$. We observe that both $u^{(1)}$ and $u^{(2)}$ have one less pair of identical consecutive letters than $u$. We may then do the same for consecutive $1$'s. By iterating such a process, we get a most frequent subword without identical consecutive letters, thus alternating between $0$ and $1$. Then we conclude by \Cref{lem:basics-4}.
\end{proof}

\begin{rmk}
  We want to highlight the importance of \Cref{prop:subword-0011} here. The main difficulty in the study of maximal subword occurrences is, in a sense, algorithmic. To the author's knowledge, we don't know whether there is a polynomial time algorithm to compute a most frequent subword of a given word, or to decide whether there is a subword that occurs at least a given number of times. However, in the case of words of the form $w = (0011)^m$, we manage to show some structure of their most frequent subwords, which then allows us to compute $\maxocc(w)$.
\end{rmk}

Let $a_{m, r} = \occ((0011)^m, (01)^r)$, and $f_{0011}(x, y) = \sum_{m,r \geq 0} a_{m, r} x^m y^r$ be their generating function. We have the following counting result.

\begin{prop} \label{prop:occ-0011}
  We have
  \[
    f_{0011}(x, y) = \frac{1-x}{(1-x)^2 - 4xy}, \quad a_{m, r} = 4^r \binom{m+r}{m-r}.
  \]
\end{prop}
\begin{proof}
  For an occurrence $P = \{p_1, \ldots, p_{2r}\}$ of $(01)^r$ in $(0011)^m$, we have two cases.
  \begin{itemize}
  \item $p_{2r} \leq 4(m - 1)$, meaning that $P$ is also an occurrence of $(01)^r$ in $(0011)^{m-1}$;
  \item $p_{2r} \in \{4m - 1, 4m\}$, meaning that the last letter $1$ of $(01)^r$ occurs at the last segment of $0011$. As the $(2r-1)$-st letter of $(01)^r$ is $0$, we have $p_{2r-1} \in \{ 4m' + 1, 4m' + 2\}$ for some $0 \leq m' \leq m - 1$. By removing both $p_{2r-1}$ and $p_{2r}$, we obtain $P'$, which is an occurrence of $(01)^{r-1}$ in $(0011)^{m'}$. To go back from $P'$ to $P$ given $m'$, we have two choices for both $p_{2r}$ and $p_{2r-1}$.
  \end{itemize}
  We thus have the recurrence for $m \geq 1$ that
  \begin{equation}
    \label{eq:rec-0011-01}
    a_{m, r} = a_{m-1, r} + \sum_{m'=0}^{m-1} 4 a_{m', r-1}.
  \end{equation}
  Subtracting \Cref{eq:rec-0011-01} for $a_{m,r}$ with that for $a_{m - 1, r}$, we have
  \[
    a_{m, r} - 2 a_{m-1, r} + a_{m-2, r} - 4 a_{m-1, r-1} = 0.
  \]
  By the standard method to convert linear recurrence with constant coefficients to rational generating function, and with the initial conditions $a_{m, 0} = 1$ and $a_{m, r} = 0$ for $r > m$, we obtain the claimed expression of $f_{0011}(x, y)$. We can then compute $a_{m, r}$ by simply extracting the coefficient of $y^r$ first, then that of $x^m$.
\end{proof}

\begin{thm} \label{thm:upper-bound}
  There is some constant $c_3$ such that, for all $n \in \naturals$, we have
  \[
    \minswentropy^{(2)}(n) \leq \frac{1}{2}\logb(1+\sqrt{2}) n - \frac{1}{2} \logb n + c_3.
  \]
\end{thm}
\begin{proof}
  For the case $n = 4m$, we have
  \begin{align*}
    \minswentropy^{(2)}(4m) \leq \swentropy((0011)^m)
    &= \max_{0 \leq r \leq m} \logb \occ((0011)^m, (01)^r) \\
    &= \frac{1}{\ln 2} \max_{0 \leq r \leq m} \ln \left(4^r \binom{m+r}{m-r}\right).
  \end{align*}
  The first equality comes from \Cref{prop:subword-0011}, and the second from \Cref{prop:occ-0011}. We take $r = \alpha m$ for some fixed $\alpha$ with $0 < \alpha < 1$. Using Stirling's approximation, we have
  \[
    \ln \left( 4^r \binom{m+r}{m-r} \right) = s(\alpha)m - \frac{1}{2} \ln m + O(1),
  \]
  where
  \[
    s(\alpha) = \alpha \ln 4 + (1+\alpha)\ln(1+\alpha) - (1-\alpha)\ln(1-\alpha) - 2\alpha\ln(2\alpha).
  \]
  The function $s(\alpha)$ is maximized at $\alpha = 2^{-1/2}$, with value $2\ln(1+\sqrt{2})$. We thus have, for some constant $c_3$, and in terms of $n = 4m$,
  \[
    \minswentropy^{(2)}(n) \leq \frac{\ln(1+\sqrt{2})}{2\ln 2} n - \frac{1}{2 \ln 2} \ln n + c_3 - \ln 4.
  \]

  For the case $n = 4m + i$ with $1 \leq i \leq 3$, let $u$ be a most frequent subword of $w = (0011)^m 010$. For an occurrence $P$ of $u$ in $w$, we take $P' = P \cap \{n-2, n-1, n\}$. Then, $j = |P'|$ can be $0, 1, 2$ or $3$. In each case, we define $u^{(j)}$ to be $u$ with the last $j$ letters removed, and there are at most 2 possibilities for $P'$. We also notice that $P \setminus P'$ is an occurrence of $u^{(j)}$ in $(0011)^m$. We thus have
  \begin{align*}
    \maxocc(w) = \occ(w, u) &= 2 \occ((0011)^m, u^{(1)}) + \occ((0011)^m, u^{(2)}) + \occ((0011)^m, u^{(3)}) \\
                            &\leq 4 \maxocc((0011)^m).
  \end{align*}
  We conclude by 
  \[
    \minswentropy^{(2)}(4m + i) \leq \swentropy^{(2)}(w) \leq \ln 4 + \swentropy^{(2)}((0011)^m) = \frac{\ln(1+\sqrt{2})}{2\ln 2} n - \frac{1}{2 \ln 2} \ln n + c_3.
  \]
  For the first inequality, we take $w'$ to be the first $(4m+i)$ letters of $w$, and it is clear that $\maxocc(w') \leq \maxocc(w)$, as each occurrence of some subword $v'$ of $w'$ is also one for $w$.
\end{proof}

The asymptotic upper bound of $\minswentropy^{(2)}(n) / n$ given by \Cref{thm:upper-bound}, thus also an upper bound of $L_2$, is $\frac{1}{2} \logb(1+\sqrt{2}) \approx 0.636\ldots$, which is much better than that in \Cref{coro:sw-bounds}. Furthermore, by regarding $(0011)^m$ as a word in a bigger alphabet, we have the following corollary, which also gives a better upper bound than that in \Cref{coro:sw-bounds}.

\begin{coro} \label{coro:upper-bound-generic}
  For all $k \geq 2$ and $n \in \naturals$, with the constant $c_3$ from \Cref{thm:upper-bound}, we have
  \[
    \minswentropy^{(k)}(n) \leq \frac{1}{2} \logb(1+\sqrt{2}) n - \frac{1}{2} \logb n + c_3.
  \]
\end{coro}

For other families of periodic words, our results follow the same procedure: first a structural result on most frequent subwords (\Cref{prop:subword-0011}), then a rational generating function for maximal subword occurrences (\Cref{prop:occ-0011}), and finally an asymptotic bound on $\minswentropy^{(2)}(n)$ by extracting the asymptotic behavior of the given rational generating function (\Cref{thm:upper-bound}).

\subsection{Subword occurrences of $(01)^m$}

\begin{prop} \label{prop:subword-01}
  For $w = (01)^m$, there is a most frequent subword $w'$ of the form $(01)^r$.
\end{prop}
\begin{proof}
  Take a most frequent subword $u$ of $w$ of length $\ell$. Suppose that there are two consecutive $0$'s in $u$, \textit{i.e.}, $u$ has the form $u = s \cdot 00 \cdot t$. We take $u' = s \cdot 010 \cdot t$. It is clear that for each occurrence of $u$ in $w$, we may associate a distinct occurrence of $u'$ by taking some extra $1$ between the occurrences of the consecutive $0$'s, as $w$ does not contain consecutive $0$'s, meaning that such a $1$ can be found. Hence, $u'$ must also be a most frequent subword of $w$. The same argument also works for consecutive $1$'s. By repeating this procedure, we can eliminate consecutive identical letters, and then we conclude by \Cref{lem:basics-4}.
\end{proof}

Let $b_{m, r} = \occ((01)^m, (01)^r)$, and $f_{01}(x, y) = \sum_{m, r \geq 0} b_{m,r} x^m y^r$ be their generating function. We have the following result.

\begin{prop} \label{prop:occ-01}
  We have
  \[
    f_{01}(x, y) = \frac{1 - x}{(1 - x)^2 - xy}, \quad b_{m, r} = \binom{m+r}{m-r}.
  \]
\end{prop}
\begin{proof}
  We first observe that, as $(01)^r$ has no consecutive identical letters, in an occurrence in $(0011)^m$, each letter of $(01)^r$ must be in different runs. By replacing the runs in $(0011)^m$ by a single letter, we obtain an occurrence of $(01)^r$ in $(01)^m$. This is a $2^{2r}$-to-$1$ bijection of occurrences of $(01)^r$ in $(0011)^m$ to those in $(01)^m$, because we have two choices for each letter of $(01)^r$ in the run of $(0011)^m$ that it occurs. We thus have $a_{m,r} = 4^r b_{m,r}$, and we conclude using \Cref{prop:occ-0011} for $b_{m, r}$ and then the standard symbolic method for $f_{01}(x, y)$.
\end{proof}

With explicit expressions of $b_{m,r}$, we can compute the asymptotic of the maximal subword occurrences of $(01)^m$.

\begin{prop} \label{prop:01}
  When $m \to \infty$, we have
  \[
    \swentropy((01)^m) = m \logb \frac{3 + \sqrt{5}}{2} - \frac{\logb m}{2} + O(1).
  \]
  The value of $r$ for the most frequent subword of the form $(01)^r$ is asymptotically $r / n = \frac{1}{\sqrt{5}}$.
\end{prop}
\begin{proof}
  From \Cref{prop:subword-01} and then \Cref{prop:occ-01}, we have
  \[
    \swentropy((01)^m) = \max_{0 \leq r \leq m} \logb \occ((01)^m, (01)^r) = \frac1{\ln 2} \max_{0 \leq r \leq m} \ln \binom{m+r}{m-r}.
  \]
  Take $r = \alpha m$. By Stirling's approximation, we have
  \[
    \ln \binom{m+r}{m-r} = s_{01}(\alpha) m - \frac{1}{2} \ln m + O(1),
  \]
  where
  \[
    s_{01}(\alpha) = (1 + \alpha) \ln(1 + \alpha) - (1 - \alpha)\ln(1 - \alpha) - 2 \alpha \ln(2 \alpha).
  \]
  The function $s_{01}(\alpha)$ is maximized at $\alpha = \frac{1}{\sqrt{5}}$, with value $\ln \frac{3 + \sqrt{5}}{2}$, which gives our claimed result.
\end{proof}

With arguments similar to that of \Cref{thm:upper-bound}, \Cref{prop:01} gives an upper bound of $L_2$ equal to $\frac{1}{2} \logb \frac{3 + \sqrt{5}}{2} \approx 0.694\ldots$, worse than that given by \Cref{thm:upper-bound}.

\subsection{Subword occurrences of $(000111)^m$}

\begin{prop} \label{prop:subword-000111}
  For $w = (000111)^m$, there is a most frequent subword $w'$ of the form $(0011)^r$.
\end{prop}
\begin{proof}
  Take a most frequent subword $u$ of $w$, and let $\ell = |u|$. By \Cref{lem:basics-4}, we can suppose that $u$ starts with $0$ and ends with $1$. We now show that, if $u$ contains one of the following factors at a certain position, then we can find another most frequent subword $u'$ of $w$ that does not contain the said factor at the same position.
  \begin{itemize}
  \item There is a run of $0$ of length $1$ in $u$, \textit{i.e.}, $u$ has the form $s \cdot 0 \cdot t$ with $s$ (resp. $t$) either empty or ending (resp. starting) with $1$. We then take $u' = s \cdot 00 \cdot t$ and show that $u'$ is also a most frequent subword of $w$. Suppose that the run of $0$ of length $1$ occurs at position $i$ of $u$. For each occurrence $P = \{p_1 < \ldots < p_\ell\}$ of $u$ in $w$, let $p_*$ be the starting position of the run of length $3$ in $w$ that contains $p_i$. It is clear that $(P \cup \{p_*, p_* + 1, p_* + 2\}) \setminus \{p_i\}$ is an occurrence of $u'$ in $w$. Furthermore, this map from occurrences of $u$ to those of $u'$ is clearly injective. We thus have $\occ(w, u') \geq \occ(w, u)$, and we conclude by the fact that $\occ(w, u)$ is maximal.
    
  \item There is a run of $0$ of length at least $3$ in $u$, \textit{i.e.}, $u$ has the form $s \cdot 000 \cdot t$. We then take $u^{(1)} = s \cdot 0 \cdot t$, $u^{(2)} = s \cdot 0100 \cdot t$ and $u^{(3)} = s \cdot 0010 \cdot t$. We show that there is some $u^{(k)}$ that is also a most frequent subword of $w$. Let $i$ be the position of the first $0$ of the three consecutive $0$'s occur in $u$. Let $P = \{p_1 < \ldots < p_\ell\}$ be an occurrence of $u$ in $w$. Let $\mathcal{P}$ (resp. $\mathcal{P}^{(k)}$ for $k = 1, 2, 3$) be the set of occurrences of $u$ (resp. $u^{(k)}$ for $k = 1, 2, 3$) in $w$. We partition $\mathcal{P}$ into $\mathcal{P}_1 \cup \mathcal{P}_2$, where $\mathcal{P}_1$ contains occurrences $P$ such that $p_{i+2} = p_{i+1} + 1 = p_i + 2$, and $\mathcal{P}_2$ contains the rest. In other words, $\mathcal{P}_1$ contains those occurrences with the three consecutive $0$'s in $u$ also consecutively in $w$, and $\mathcal{P}_2$ contains the others. Given $P \in \mathcal{P}_1$, we see that $(P \setminus \{p_i, p_{i+1}, p_{i+2}\}) \cup \{p_{i+j}\}$ is in $\mathcal{P}^{(1)}$ for all $0 \leq j \leq 2$, and this is a $1$-to-$3$ injection from $\mathcal{P}_1$ to $\mathcal{P}^{(1)}$. We thus have $|\mathcal{P}^{(1)}| \geq 3 |\mathcal{P}_1|$. Now, given $P \in \mathcal{P}_2$, we know that $p_i, p_{i+1}, p_{i+2}$ are not in the same run in $w$, so we can take the smallest $p_*$ such that $p_i < p_* < p_{i+2}$ and $w_{p_*} = 1$. As $p_*$ is the smallest such index, it is also the beginning of a run of $1$ of length $3$ in $w$. Hence, $P \cup \{p_* + j\}$ is an occurrence of $u^{(2)}$ or $u^{(3)}$ for all $0 \leq j \leq 2$. This is a $1$-to-$3$ injection from $\mathcal{P}_2$ to $\mathcal{P}^{(2)} \cup \mathcal{P}^{(3)}$, meaning that $|\mathcal{P}^{(2)}| + |\mathcal{P}^{(3)}|  \geq 3 |\mathcal{P}_2|$. We thus have
    \[
      \sum_{k=1}^3 \occ(w, u^{(k)}) = \sum_{k=1}^3 |\mathcal{P}^{(k)}| \geq 3 |\mathcal{P}_1| + 3 |\mathcal{P}_2| = 3 \occ(w, u).
    \]
    We then conclude by the fact that $u$ is a most frequent subword of $w$.
  \end{itemize}

  We notice that the reasoning above also works for runs of $1$ of length $1$ or at least $3$, and the new most frequent subword has the same starting and ending letter as $u$. Starting from $u$, we first apply successively the replacement that removes runs of length 1 until no such run is left. We denote by $u'$ the resulting word, which remains a most frequent subword of $w$. Then we apply to $u'$ the replacement for runs of length at least 3, and when a run of length 1 occurs, we perform a replacement to remove it. Such a combined operation does not create new runs of length 1, and the resulting word either has strictly more runs, or has the same number of runs but being strictly shorter. By iteratively applying this combined operation, we obtain a sequence of most frequent subwords that is increasing in the number of runs. As the number of runs is bounded by $|w|$, it becomes stationary in the sequence after a certain point, after which the length of the most frequent subwords obtained decreases. The procedure thus eventually stops, and the word obtained at the end, denoted by $w'$, is a most frequent subword of $w$ that starts with $0$, ends with $1$, and does not contain runs of length $1$ or at least $3$. It is thus of the form $(0011)^r$.
\end{proof}

Let $c_{m, r} = \occ((000111)^m, (0011)^r$, and $f_{000111}(x, y) = \sum_{m, r \geq 0} c_{m,r} x^m y^r$ be their generating function. We have the following expression of $f_{000111}$.

\begin{prop} \label{prop:occ-000111}
  We have
  \[
    f_{000111}(x, y) = \frac{(1-x)^3}{(1-x)^4 - 9x (1+2x)^2 y}.
  \]
\end{prop}
\begin{proof}
  For an occurrence $P = \{p_1, \ldots, p_{4r}\}$ of $w' = (0011)^r$ in $w = (000111)^m$, we have the following cases.
  \begin{itemize}
  \item $p_{4r} \leq 6(m - 1)$, meaning that $P$ is also an occurrence of $w'$ in $(000111)^{m-1}$.
  \item $p_{4r} \in \{ 6m, 6m - 1, 6m - 2\}$, meaning that the last letter of $w'$ occurs in the last run of $w$. As the $(4r-3)$-th letter of $w'$ is $0$, we have $p_{4r-3} \in \{ 4m' + 1, 4m' + 2, 4m' + 3 \}$ for some $m'$ with $0 \leq m' \leq m - 1$. Then, $P' = \{p_1, \ldots, p_{4r-4}\}$ is an occurrence of $(0011)^{r-1}$ in $(000111)^{m'}$. Given $P'$, we count the number of possibilities to reconstruct $P$. For $p_{4r-3}$ and $p_{4r-2}$, which are all occurrences of $0$, we have the following cases:
    \begin{itemize}
    \item $p_{4r-2} \in \{ 4m' + 1, 4m' + 2, 4m' + 3 \}$, meaning that the two positions are in the same run of length $3$, leading to $3$ possibilities. For the $p_{4r-1}$ and $p_{4r}$, either they are in the same run, leading to $3$ possibilities in total, or they are in different runs, leading to $3$ choices for $p_{4r}$ and $3(m - m' - 1)$ for $p_{4r-1}$, as there are $m - m'$ runs of $1$ after $p_{4r-2}$, but $p_{4r-1}$ cannot be in the one containing $p_{4r}$. Therefore, there are totally $3(3m - 3m' - 2)$ choices for $p_{4r-1}$ and $p_{4r}$.
    \item $p_{4r-2} \in \{ 4m'' + 1, 4m'' + 2, 4m'' + 3 \}$ for some $m''$ such that $m' < m'' < m$, meaning that the two positions are in different runs of length $3$. In this case, using the same argument as the point above, we know that there are totally $3(3m - 3m'' - 2)$ choices for $p_{4r-1}$ and $p_{4r}$. For $p_{4r-3}$ and $p_{4r-2}$, we have $3$ choices for each.
    \end{itemize}
  \end{itemize}
  The case analysis above thus leads to the following recurrence:
  \begin{align*}
    c_{m,r} &= c_{m-1, r} + 9 \sum_{m' = 0}^{m-1} (3m - 3m' - 2) c_{m', r-1} + 9 \sum_{m' = 0}^{m - 1} \sum_{m'' = m' + 1}^{m - 1} 3(3m - 3m'' - 2) c_{m', r-1} \\
            &= c_{m-1, r} + 9 \sum_{m' = 0}^{m-1} (3m - 3m' - 2) c_{m', r-1} + 9 \sum_{m' = 0}^{m - 1} \left( \frac{9}{2}(m - m')^2 - \frac{21}{2}(m - m') + 6 \right) c_{m', r-1} \\
            &= c_{m-1, r} + 9 \sum_{m' = 0}^{m - 1} \left( \frac{9}{2}(m - m')^2 - \frac{15}{2}(m - m') + 4 \right) c_{m', r-1}.
  \end{align*}
  Take the recurrence above minus the same recurrence for $c_{m - 1, r}$, we have
  \[
    c_{m, r} = 2 c_{m-1, r} - c_{m-2, r} + 9 \sum_{m'=0}^{m-2} (9(m - m') - 12) c_{m', r-1} + 9 c_{m-1, r-1}.
  \]
  Again, take the equation above and subtract it by itself with $m$ replaced by $m - 1$, we have
  \[
    c_{m, r} = 3 c_{m-1, r} - 3 c_{m-2, r} + c_{m-3, r} + 81 \sum_{m'=0}^{m-3} c_{m', r-1} + 45 c_{m-2, r-1} + 9 c_{m-1, r-1}.
  \]
  Doing the same for one last time, we have
  \[
    c_{m, r} = 4 c_{m-1, r} - 6 c_{m-2, r} + 4 c_{m-3, r} - c_{m-4, r} - 36 c_{m-3, r-1} + 36 c_{m-2, r-1} + 9 c_{m-1, r-1}.
  \]
  The recurrence above is valid for all $m \geq 1$ if we take the convention that $c_{m, r} = 0$ for all $m < 0$. With the boundary condition $c_{m, 0} = 1$ for all $m \geq 0$, by the standard method for converting linear recurrences into rational generating function, we have our claim.
\end{proof}

To obtain the asymptotic maximal occurrences of $(000111)^m$, we will use the following result on saddle-point estimates of large powers \cite[Theorem~VIII.8]{flajolet}. Note that we may also use the ACSV approach \cite{melczer,mishna}, which can even be automated \cite{acsv-sage}.

\begin{thm}[Special version of {\cite[Theorem~VIII.8]{flajolet}}] \label{thm:large-power}
  Suppose that $A(x)$, $B(x)$ are power series in $x$ with non-negative coefficients convergent in a neighborhood of $x = 0$, while satisfying the following conditions:
  \begin{itemize}
  \item $B(0) \neq 0$, and $B(x)$ is aperiodic, \textit{i.e.}, $B(x)$ is not a function of the form $\beta(x^p)$ for some power series $\beta$ also convergent in a neighborhood of $x = 0$;
  \item The radius of convergence of $B(x)$, denoted by $R$, is not larger than that of $A(x)$.
  \end{itemize}
  Then, for $N = \lambda n$ with some value of $\lambda > 0$, when $n \to \infty$, we have
  \begin{equation} \label{eq:large-power}
    [x^N] A(x)B(x)^n = A(\zeta) \frac{B(\zeta)^n}{\zeta^{N+1}\sqrt{2\pi n \xi}} (1 + o(1)).
  \end{equation}
  Here, $\zeta$ is the unique positive root of $\zeta B'(\zeta) = \lambda B(\zeta)$, and $\xi$ an explicit constant. The estimate \eqref{eq:large-power} holds uniformly for $\lambda$ in any compact in $(0, T)$, with $T$ the limit of $xB'(x)/B(x)$ for $x \to R^-$.
\end{thm}

\begin{prop} \label{prop:000111}
  When $m \to \infty$, we have
  \[
    \swentropy((000111)^m) = m \gamma - \frac{\logb m}{2} + O(1),
  \]
  where
  \[
    \gamma = \alpha \logb 9 + 2\alpha \logb \frac{1+2\zeta}{(1-\zeta)^2} - (1 - \alpha) \logb \zeta \approx 3.9215913\ldots.
  \]
  The constant $\alpha \approx 0.6597177\ldots$ is the unique positive solution of
  \begin{equation} \label{eq:alpha-sol}
    457\alpha^4 - 246\alpha^2 + 72\alpha - 27 = 0.
  \end{equation}
  The constant $\zeta$ is given as an expression of $\alpha$:
  \begin{equation} \label{eq:zeta-sol}
    \zeta = \frac{1 - 9 \alpha + \sqrt{73 \alpha^2 - 18 \alpha + 9}}{4 + 4 \alpha}.
  \end{equation}
\end{prop}
\begin{proof}
  We set $\alpha = r / m$, and we want to find the value of $\alpha$ that maximizes $c_{m,r}$ asymptotically. By \Cref{conj:half-length}, we have
  \[
    c_{m, r} = [x^m y^r] \frac{(1-x)^3}{(1-x)^4 - 9x (1+2x)^2 y} = [x^{m-r}] \frac{9^r}{1-x} \left(\frac{1+2x}{(1-x)^2}\right)^{2r}.
  \]
  We apply \Cref{thm:large-power} with $A(x) = (1-x)^{-1}$, $B(x) = (1+2x)(1-x)^{-2}$, $n = 2r$ and $N = m - r$. We check that $A(x)$ and $B(x)$ satisfy the given conditions, and we have $\lambda = N / n = (\alpha^{-1} - 1) / 2$. When $m \to \infty$ with fixed $\alpha$, by substituting $n, N, r$ with expressions in $m$ and $\alpha$, we have
  \[
    \logb c_{m, r} = \alpha m \logb 9 + 2 \alpha m \logb B(\zeta) - (1 - \alpha) m \logb \zeta - \frac{\logb m}{2} + O(1).
  \]
  Here, $\zeta$ is the unique positive solution of $2x (x+2) - \lambda (1+2x) (1-x) = 0$, whose expression in $\alpha$ is exactly that given in \eqref{eq:zeta-sol}. We note that $xB'(x)/B(x) = \frac{2x(x+2)}{(1+2x)(1-x)}$ tends to infinity when $x \to 1^{-}$, which means that the estimation is valid uniformly for $\lambda$ in any compact in $(0, +\infty)$, which means for $\alpha$ in any compact in $(0, 1)$.
  
  Now, we try to find the value of $\alpha$ that maximizes $\logb c_{m,r}$ asymptotically, which means maximizing the following function:
  \[
    C(\alpha) = \alpha \logb 9 + 2 \alpha \logb B(\zeta) - (1-\alpha) \logb \zeta.
  \]
  By substituting the expression of $\zeta$ in \eqref{eq:zeta-sol} and differentiating, we see that the value of $\alpha$ that maximizes $\logb c_{m,r}$ asymptotically is the unique positive solution of the equation \eqref{eq:alpha-sol}. An exact expression of the solution in radicals can be obtained by solving the quartic equation above. However, we choose not to display it here, as it is quite complicated but not so useful for our purpose. We thus have our result.
\end{proof}

Again, with arguments similar to that of \Cref{thm:upper-bound}, \Cref{prop:000111} gives an upper bound of $L_2$ approximately $0.6536\ldots$, again worse that that given by \Cref{thm:upper-bound}.

\subsection{Subword occurrences of periodic words in another periodic word}

To find better upper bounds, it is natural to try to look at other families of periodic words. This is encouraged by the following theorem.

\begin{thm} \label{thm:periodic-family}
  For any words $w, v$ in an alphabet $\mathcal{A}$ of size $k$, the generating function $f_{w, v}(x, y) = \sum_{m, r \geq 0} \occ(w^m, v^r) x^m y^r$ is rational in $x, y$.
\end{thm}
\begin{proof}
  We define $a_{w,v}^{s,t}(m)$ with $1 \leq s, t \leq |w|$ to be the number of occurrences $P = \{p_1, \ldots, p_{|v|} \}$ of $v$ in $w^m$ such that $p_1 = s$ and $p_{|v|} = (m-1)|w| + t$. In other words, $a_{w,v}^{s,t}(m)$ counts the occurrences of $v$ in $w^m$ such that the first (resp. last) letter of $v$ occurs in the first (resp. last) copy of $w$ at position $s$ (resp. $t$). Let $g_{w, v}^{s, t}(x) = \sum_{m \geq 1} a_{w,v}^{s,t}(m) x^{m-1}$. Note the extra $-1$ in the exponent of $x$ in $g_{w,v}^{s,t}(x)$. We first show that $g_{w, v}^{s, t}(x)$ is rational. For an occurrence $P$ of $v$ in $w^m$, some consecutive letters may occur in the same copy of $w$. We say that such letters form a cluster, and we denote by $\sigma$ the integer composition of the number of letters in each cluster from left to right. We denote by $\ell(\sigma)$ the length of $\sigma$, which is also the number of clusters. We denote the clusters by $v_\sigma^{(1)}, \ldots v_\sigma^{(\ell(\sigma))}$, and it is clear that they are obtained by cutting $v$ into pieces whose lengths are the parts of $\sigma$. We then have
  \begin{align*}
    &g_{w, v}^{s, t}(x) = a_{w,v}^{s,t}(1)  \\
    &+ \sum_{m \geq 2} \sum_{\sigma \vDash |v|} x^{m-1} \binom{m - 2}{\ell(\sigma) - 2} \left( \sum_{t' = s}^{|w|} a_{w,v_\sigma^{(1)}}^{s,t'}(1) \right) \left( \sum_{s' = 1}^{t} a_{w,v_\sigma^{(\ell(\sigma))}}^{s',t}(1) \right) \prod_{i=2}^{\ell(\sigma) - 1} \occ(w, v_\sigma^{(i)}).
  \end{align*}
  Here, $\sigma \vDash |v|$ means that we go over all integer compositions of $|v|$. The first term is for $m=1$. For the second term, we simply count all possibilities of how clusters of $v$ appear in $w^m$ with $m \geq 2$ while fixing the first and the last cluster. We observe that each $a_{w,v_{\sigma^{(i)}}}^{s', t'}(1)$ for any $s', t', i$ is a constant, and the same holds for $\occ(w, v_{\sigma}^{(i)})$. By exchanging the two summations, and observing that $\sum_{m \geq 2} \binom{m-2}{d-2} x^{m-1} = x^{d-1} (1-x)^{-(d-1)}$, we see that $g_{w,v}^{s,t}(x)$ is rational in $x$ with $(1 - x)^{|v|-1}$ as denominator, as $\ell(\sigma) \leq |v|$ for $\sigma \vDash |v|$.

  Now, for $1 \leq t \leq |w|$, we define $f_{w,v}^{(t)}(x, y) = \sum_{m \geq 1} \sum_{r \geq 1} b_{w, v}^{(t)}(m, r) x^{m-1} y^r$ with $b_{w, v}^{(t)}(m, r)$ counting the number of occurrences $P = \{ p_1, \ldots p_{|v|r} \}$ of $v^r$ in $w^m$ such that $p_{|v|r} = (m-1)|w| + t$. Again, we note the extra $-1$ in the exponent of $x$. We see that $b_{w, v}^{t}(m, r)$ is defined similarly as $a_{w, v}^{s,t}(m)$, except that we consider subwords of the form $v^r$, and we do not fix the position of the first letter of $v^r$ in $w^m$. We thus have $b_{w, v}^{(t)}(m, 1) = \sum_{m'=1}^m \sum_{s=1}^{|w|} a_{w, v}^{s, t}(m')$. Now, let $P$ be an occurrence of $v^r$ in $w^m$ counted by $b_{w, v}^{(t)}(m, r)$. By considering the copies of $w$ spanned by the last copy of $v$, we have
  \begin{align*}
    f_{w,v}^{(t)}(x, y) = \frac{y}{1-x} \sum_{s=1}^{|w|} g_{w,v}^{s,t}(x)
    &+ \frac{xy}{1-x} \left( \sum_{t'=1}^{|w|} f_{w,v}^{(t')}(x, y) \right) \left( \sum_{s=1}^{|w|} g_{w,v}^{s,t}(x) \right) \\
    &+ y \sum_{t'=1}^{|w| - 1} \left( f_{w,v}^{(t')}(x, y) \sum_{s=t'+1}^{|w|} g_{w,v}^{s,t}(x)  \right).
  \end{align*}
  Here, the first term is for $r = 1$, with the factor $(1-x)^{-1}$ for the copies of $w$ before positions in which letters of $v$ occur, and the rest is for $r \geq 2$. There are two cases: either letters in the $r$-th and the $(r-1)$-st copies of $v$ do not occur in the same copy of $w$ in $w^m$, or they do. The first case is counted by the second term above, with the factor $(1-x)^{-1}$ for copies of $w$ between the occurrences of the two last copies of $v$ in $w^m$. The second case is accounted by the third term above, where we have the constraint that the last letter of the $(r-1)$-st copy of $v$ occurs before the first letter of the $r$-th copy in the same copy of $w$.

  Let $\mathbf{f} = {}^t (f_{w,v}^{(1)}, \ldots, f_{w,v}^{(|w|)})$. The equation above can be seen as $\mathbf{Af} = \mathbf{b}$ for some matrix $\mathbf{A} = (A_{i, j})_{1 \leq i, j \leq |w|}$ and some row vector $\mathbf{b}$, both with coefficients that are linear in $y$ and rational in $x$, and with only powers of $(1-x)$ as denominators. We also observe that $A_{i, i}$ is of the form $1 + R(x)y$ with $R(x)$ rational in $x$, while $A_{i, j}$ for $i \neq j$ is of the form $R(x)y$. Hence, $\mathbf{A}$ is non-singular, and $f_{w,v}^{(i)}$ is rational in $x, y$ for all $1 \leq i \leq |w|$. We conclude by observing that $f_{w,v}(x, y) = \frac{1}{1 - x} + \frac{x}{1 - x} \sum_{t=1}^{|w|} f_{w,v}^{(t)}(x, y)$, with the $\frac{1}{1 - x}$ taking care of the case $m=0$ and the case $r = 0$, the factor $x$ to make up for the definition of $f_{w,v}^{(t)}$, and the factor $(1-x)^{-1}$ for the copies of $w$ after the last cluster of $v^r$.
\end{proof}

In principle, for any word $w$ and $v$, we can first compute $f_{w, v}(x, y)$ effectively as in the proof of \Cref{thm:periodic-family}, then use analytic combinatorics in several variables \cite{melczer,mishna} to compute the asymptotically maximal value of $\occ(w^m, v^r)$ for fixed $m$. Although the computation of $f_{w, v}(x, y)$ would be tedious, it can be automated, and an implementation of \Cref{thm:periodic-family} in Sagemath is provided in \cite{subword-repo}. The only problem is that, for $w$ in general, we do not have results like \Cref{prop:subword-0011} for the structure of most frequent subwords of $w^m$, meaning that $\maxocc(w^m)$ is not necessarily achieved for subwords of the form $v^r$.

\section{Experimental results and open questions} \label{sec:open}

Generally, the ``minimum of maximums'' structure in the definition of $\minswentropy^{(k)}(n)$ makes estimates difficult. Hence, not a lot is known about $\minswentropy^{(k)}(n)$. We thus performed some computational experiments to compute the binary case $\minswentropy^{(2)}(n)$ for small $n$, and also for computing most frequent subwords of some given words. Such experiments also inspired our theoretical results presented earlier. For instance, the bounds in \Cref{sec:binary} were first observed experimentally before proven rigorously.

For our experiments, we have developed a program in C that performs various computations concerning maximal subword occurrences \cite{subword-repo}. Our program is mostly specialized for computing $\minswentropy^{(2)}(n)$, but some functions may be adapted for other computations concerning subword occurrences.

To find the binary words of a given length $n$, there are several levels of computations:
\begin{enumerate}
\item For binary words $w, w'$, compute the number of occurrences $\occ(w, w')$ of $w'$ in $w$;
\item For a binary word $w$, find its most frequent subwords;
\item Find the binary words $w$ of length $n$ with the minimal value of $\maxocc(w)$.
\end{enumerate}

For the first level, we represent $w$ and $w'$ by a tuple of the lengths of their runs that we call the \emph{run-length tuple}. For instance, the word $0000110111001$ is represented by $(4, 2, 1, 3, 2, 1)$ By the symmetry of letters $0$ and $1$, and by \Cref{lem:basics-4}, we may assume that both $w$ and $w'$ starts with the letter $0$. The length of the runs thus uniquely determine such words. Furthermore, \Cref{lem:basics-4} also implies that we only need to consider $w'$ whose parity of number of runs is the same of that of $w$. Suppose that $(\ell_1, \ldots, \ell_m)$ (resp. $(\ell'_1, \ldots, \ell'_r)$) the run-length tuple of $w$ (resp. $w'$). To accelerate the computation of $\occ(w, w')$, we use a divide-and-conquer strategy. Take $r_* = \lfloor (r+1)/2 \rfloor$. We first consider the possible occurrences of the letters in the $r_*$-th run of $w'$. Suppose that the first and the last letter of this run of $w'$ occur in the $k$-th and $k'$-th runs of $w$ with $k \leq k'$, where $k$ and $k'$ has the same parity as $r_*$. Let $u$, $v$, $u'$, $v'$ be the word represented by run-length tuples $(\ell_1, \ldots, \ell_{k-1})$, $(\ell_{k'+1}, \ldots, \ell_m)$, $(\ell'_1, \ldots, \ell'_{r_*-1})$ and $(\ell'_{r_*+1}, \ldots, \ell'_r)$ respectively. We take $L = \sum_{j = 0}^{(k'-k)/2} \ell_{k+2j}$, and the number of the subset of occurrences of $w'$ in $w$ we are considering is given by
\[
  \occ(u, u') \occ(v, v') \left[ \binom{L}{\ell_{r_*}} - \binom{L - \ell_k}{\ell_{r_*}} - \binom{L - \ell_{k'}}{\ell_{r_*}} + \binom{L - \ell_k - \ell_{k'}}{\ell_{r_*}} \right].
\]
Here, we use the inclusion-exclusion principle to ensure the positions of the first and the last letter of the $r_*$-th run of $w'$ in $w$. The case $k = k'$ needs a special treatment. The values of $\occ(u, u')$ and $\occ(v, v')$ can be computed recursively. To obtain $\occ(w, w')$, it suffices to sum over all possible $k$ and $k'$. Such a divide-and-conquer strategy is more efficient than scanning runs from left to right, as the sub-problems to be computed have much smaller size, and incoherence can be detected much faster. To avoid recomputing the same sub-problems, we use the memoization technique, storing results of words whose number of runs is below a given threshold, and reuse them later. This technique gives considerable speedup, as we will compute the number of occurrences of a lot of (similar) subword in all the words of a given length, meaning that almost all possible combinations of $w$ and $w'$ with few runs are reused many times.

For the second level, we need to generate subwords of all lengths and compute their occurrences in the given word $w$. However, some lengths give subwords that may occur much more often in $w$ than others. For instance, subwords of length $1$ may occur in $w$ for at most $|w|$ times, which is small. To obtain faster the most frequent subwords, a first heuristic is to check the lengths closer to $|w|/2$ first, but we may use other hints given by the third level. It may seem that the order we check the subwords is not important, as we need to go through all subwords to obtain the most frequent ones. However, as we will see in the following, this exhaustion is not needed for the third level. This realization gives a significant speedup, as there are exponentially many subwords to check if we do not skip some of them.

For the third level, we go over all words of length $n$ as binary representation of integers from $0$ to $2^{n-1}$, as we may check only words starting with the letter $0$. We then make use of the ``minimum of maximums'' structure of our problem. Suppose that we have computed $\maxocc(w)$ for some word $w$ of length $n$. When looking at another word $u$ of the same length, if we can find a subword $u'$ such that $\occ(u, u') > \maxocc(w)$, then we do not need to look at other subwords of $u$, as we already have $\maxocc(u) \geq \occ(u, u') > \maxocc(w)$, meaning that $u$ cannot be a word realizing $\minswentropy^{(2)}(n)$. Hence, before performing the full second level computation for the word $u$, we can try to find a good subword $u'$ that occurs frequent enough in $u$ to rule it out. As consecutive words are similar, heuristically, most frequent subwords make good candidates to rule out the current word. Hence, at each time when we need to perform the second level computation on a word $u$, we save one of the most frequent subwords $u'$. Then, when we check a subsequent word $v$, we first check whether $u'$ or some of its variants occurs enough times in $v$ to rule it out. This heuristic gives significant speedup to the whole computation. Even if we need to go into the second level computation for $v$, we still can use the length of $u'$ as a hint by first searching for subwords $v'$ with length close to that of $u'$. We also update the record of $\maxocc(w)$ when we find a word $v$ with $\maxocc(v) < \maxocc(w)$.

To further speed up the computation, we perform a meta-heuristic search to find words $w$ with relatively small $\maxocc(w)$ in order to eliminate more words at the beginning of the third level computation. The search algorithm we use is a home-brew combination of exhaustive local search and stochastic long jumps. More precisely, we perform at each time an exhaustive local search from the current best result to try to improve $\maxocc(w)$. When no improvement is seen, we randomly flip some bits before the exhaustive local search. After a given number of such attempts of flipping, if there is still no improvement, we increase the rate of flipping each bit, until reaching a given upper bound and terminate. If an improvement is found, we reset the rate of bit-flipping. This strategy is inspired by self-adaptive parameter control in meta-heuristics (\textit{c.f.} \cite{metaheuristic}). Another way to find words $w$ of length $n$ with small $\maxocc(w)$ is to take the words $w'$ of length $n - 1$ achieving minimal subword entropy, insert a new letter in all possible ways, then find the way that produces a word $w$ with the lowest subword entropy. This method, suggested by Jim Fill in personal communication, is more efficient than the meta-heuristic search, while giving results of similar quality.

We also harness the power of modern processors that have multiple cores by using several threads to go through words in parallel. Concurrent data structure is needed, especially in the memoization for accelerating the computation of subword occurrences.

We now present experimental results obtained using our program. We denote by $\overline{w}$ the word obtained from $w$ by switching $0$ and $1$, and $\overleftarrow{w}$ the reverse of $w$. By symmetry between the two letters, we have the following simple observation.

\begin{lem} \label{lem:switch}
  For any $w \in \{0, 1\}^n$ with $n \geq 0$, we have $\maxocc(w) = \maxocc(\overline{w}) = \maxocc(\overleftarrow{w})$.
\end{lem}

\begin{table}
  \centering \small
  \begin{tabular}[t]{lllrll}
    \toprule
    $n$ & Words $w$ with lowest $\swentropy^{(2)}(w)$ & $\maxocc(w)$ & $\swentropy^{(2)}(w)$ & $\swentropy^{(2)}(w)/n$ & \#runs \\
    \midrule
    $1$ & 0 & $1$ & $0$ & $0$ & $1$ \\
    $2$ & 01 & $1$ & $0$ & $0$ & $2$ \\
    $3$ & 001 & $2$ & $1$ & $0.333$ & $2$ \\
    & 010 & & & & $3$ \\
    $4$ & 0110 & $2$ & $1$ & $0.25$ & $3$ \\
    $5$ & 01110 & $3$ & $1.585$ & $0.317$ & $3$ \\
    $6$ & 011001 & $5$ & $2.322$ & $0.387$ & $4$ \\
    $7$ & 0110001 & $6$ & $2.585$ & $0.369$ & $4$ \\
    $8$ & 01110001 & $9$ & $3.170$ & $0.396$ & $4$ \\
    $9$ & 011000110 & $16$ & $4$ & $0.444$ & $5$ \\
    $10$ & 0110001110 & $22$ & $4.459$ & $0.446$ & $5$ \\
    $11$ & 01110001110 & $33$ & $5.044$ & $0.459$ & $5$ \\
    $12$ & 011000111001 & $52$ & $5.700$ & $0.475$ & $6$ \\
    $13$ & 0111001001110 & $72$ & $6.170$ & $0.475$ & $7$ \\
    $14$ & 01100010111001 & $108$ & $6.755$ & $0.482$ & $8$ \\
    $15$ & 011000101110001 & $162$ & $7.340$ & $0.489$ & $8$ \\
    $16$ & 0111000101110001 & $252$ & $7.977$ & $0.499$ & $8$ \\
    $17$ & 01100011111000110 & $390$ & $8.607$ & $0.506$ & $7$ \\
    $18$ & 011100100101110001 & $588$ & $9.200$ & $0.511$ & $10$ \\
    $19$ & 0110001011101000110 & $900$ & $9.814$ & $0.517$ & $11$ \\
      & 0110001110110001110 & & & & \\
    $20$ & 01110001011011000110 & $1320$ & $10.366$ & $0.518$ & $11$ \\
    $21$ & 011100011011010001110 & $2049$ & $11.000$ & $0.524$ & $11$ \\
    $22$ & 0110001110101000111001 & $2958$ & $11.530$ & $0.524$ & $12$ \\
    $23$ & 01110001011011010001110 & $4473$ & $12.127$ & $0.527$ & $13$ \\
    $24$ & 011000111010101000111001 & $6979$ & $12.769$ & $0.532$ & $14$ \\
    $25$ & 0111000101101101000111001 & $10602$ & $13.372$ & $0.535$ & $14$ \\
    $26$ & 01110001011011001000111001 & $15962$ & $13.962$ & $0.537$ & $14$ \\
    $27$ & 011100010101110101000111001 & $24150$ & $14.560$ & $0.539$ & $16$ \\
    $28$ & 0110001111010010010111000110 & $36450$ & $15.154$ & $0.541$ & $15$ \\
        & 0111000101110101000101110001 & & & & $16$ \\
    $29$ & 01100011101010001010111000110 & $53671$ & $15.712$ & $0.542$ & $17$ \\
    $30$ & 011000111001100010101111000110 & $83862$ & $16.356$ & $0.545$ & $15$ \\
    $31$ & 0110001110101000101011110001110 & $127998$ & $16.966$ & $0.547$ & $17$ \\
    $32$ & 01100011101010001010111010001110 & $189131$ & $17.529$ & $0.548$ & $19$ \\
    $33$ & 011000111101010001011011010001110 & $288900$ & $18.140$ & $0.550$ & $19$ \\
    $34$ & 0110001110101000101011101001001110 & $442386$ & $18.755$ & $0.552$ & $21$ \\
    $35$ & 01110001011011001000110111001001110 & $681966$ & $19.379$ & $0.554$ & $19$ \\
    $36$ & 011100010111010100010110111001001110 & $1047330$ & $19.998$ & $0.556$ & $21$ \\
    $37$ & 0111000101101011000011011011010001110 & $1581150$ & $20.593$ & $0.557$ & $21$ \\
    $38$ & 01110001011011011000100111011001001110 & $2387054$ & $21.187$ & $0.558$ & $21$ \\
    $39$ & 011000110110010011101100010010111000110 & $3626580$ & $21.790$ & $0.559$ & $21$ \\
    $40$ & 0110001110101000101011101010001110010110 & $5500610$ & $22.391$ & $0.560$ & $25$ \\
    \bottomrule
  \end{tabular}
  \caption{Binary words achieving minimal subword entropy of length from $1$ to $40$. In each equivalent class defined by the symmetries in \Cref{lem:switch}, only one representative is given. Numerical values are rounded to three digits after the decimal point when needed.}
  \label{tab:exhaustive}
\end{table}

\begin{table}
  \centering \small
  \begin{tabular}[t]{llll}
    \toprule
    $n$ & Words $w$ with lowest $\swentropy^{(2)}(w)$ & Most frequent subwords $w'$ & $|w'|/|w|$ \\
    \midrule
    $1$ & 0 & 0 & 1 \\
    $2$ & 01 & 01 & 1 \\
    $3$ & 001 & 0, 01 & 0.33--0.67 \\
    & 010 & 0 & 0.33 \\
    $4$ & 0110 & 0 & 0.25 \\
    $5$ & 01110 & 010, 0110 & 0.6--0.8 \\
    $6$ & 011001 & 01 & 0.33 \\
    $7$ & 0110001 & 01, 001, 0101, 01001 & 0.29--0.71 \\
    $8$ & 01110001 & 0101, $01001^\bullet$, 011001 & 0.5--0.75 \\
    $9$ & 011000110 & 010 & 0.33 \\
    $10$ & 0110001110 & 0110 & 0.4 \\
    $11$ & 01110001110 & 0110 & 0.33 \\
    $12$ & 011000111001 & 0101 & 0.33 \\
    $13$ & 0111001001110 & 01010, 010010, $010110^*$ & 0.38--0.62 \\
     & & $0100110^*$, 0110110, 01100110 & \\
    $14$ & 01100010111001 & 010101, $0100101^\bullet$, 01001101 & 0.43--0.57 \\
    $15$ & 011000101110001 & 0101101, 01001101 & 0.47--0.6\\
     & & 01011001, 010011001 & \\
    $16$ & 0111000101110001 & 011001 & 0.38 \\
    $17$ & 01100011111000110 & $0100110^*$, $0101110^*$ & 0.41 \\
    $18$ & 011100100101110001 & 01001101, 010011001 & 0.44--0.5 \\
    $19$ & 0110001011101000110 & 01011010, $010011010^*$, 0100110010 & 0.42--0.53 \\
      & 0110001110110001110 & 010110110 & \\
    $20$ & 01110001011011000110 & 010011010 & 0.45 \\
    $21$ & 011100011011010001110 & 01100110 & 0.38 \\
    $22$ & 0110001110101000111001 & $010011001^\bullet$ & 0.41 \\
    $23$ & 01110001011011010001110 & 0100110010 & 0.43 \\
    $24$ & 011000111010101000111001 & $010101001^\bullet$ & 0.38 \\
    $25$ & 0111000101101101000111001 & 0110011001 & 0.4 \\
    $26$ & 01110001011011001000111001 & 01001100101 & 0.42 \\
    $27$ & 011100010101110101000111001 & 0100110101 & 0.37 \\
    $28$ & 0110001111010010010111000110 & 01100110010 & 0.39--0.43\\
        & 0111000101110101000101110001 & 010011001101 & \\
    $29$ & 01100011101010001010111000110 & 0101001010 & 0.34 \\
    $30$ & 011000111001100010101111000110 & 010110010110 & 0.4 \\
    $31$ & 0110001110101000101011110001110 & 0110100110110 & 0.42 \\
    $32$ & 01100011101010001010111010001110 & 0101000110010 & 0.41 \\
    $33$ & 011000111101010001011011010001110 & 0101101010110 & 0.39 \\
    $34$ & 0110001110101000101011101001001110 & 0110101010110 & 0.38 \\
    $35$ & 01110001011011001000110111001001110 & 01100101100110 & 0.4 \\
    $36$ & 011100010111010100010110111001001110 & 0110011010110110 & 0.44 \\
    $37$ & 0111000101101011000011011011010001110 & 0100101001110010 & 0.43 \\
    $38$ & 01110001011011011000100111011001001110 & 0100111001110110 & 0.42 \\
    $39$ & 011000110110010011101100010010111000110 & 01010110001010 & 0.36 \\
    $40$ & 0110001110101000101011101010001110010110 & 010100101100110 & 0.38 \\
    \bottomrule
  \end{tabular}
  \caption{Binary words achieving minimal subword entropy of length from $1$ to $40$ and some of their most frequent subwords, modulo equivalent classes defined by the symmetries in \Cref{lem:switch}. We only show subwords that start and end with the same letters at the word $w$. The symbol $*$ (resp. $\bullet$) is for subwords that are not palindromic (resp. anti-palindromic) while the word $w$ is. Only one representative of subwords is shown in these two cases. Numerical values of length ratio are rounded to two digits.}
  \label{tab:exhaustive-subwords}
\end{table}

We now give the words achieving minimal subword entropy of length up to $40$ in \Cref{tab:exhaustive} with information on subword entropy and runs, up to the symmetries in \Cref{lem:switch}. We then give the information about some of their most frequent subwords in \Cref{tab:exhaustive-subwords}. The case $n = 40$ took around $6$ days on a local computation server with all its $64$ threads on $32$ cores. There are several observations we can draw from \Cref{tab:exhaustive} and \Cref{tab:exhaustive-subwords}, but few is without exception.

\begin{itemize}
\item The words of length $n$ achieving $\minswentropy^{(2)}(n)$ are palindromic, \textit{i.e.}, $w = \overleftarrow{w}$, or anti-palindromic, \textit{i.e.}, $\overline{w} = \overleftarrow{w}$, for many values of $n$, such as 1, 2, 4, 5, 6, 8, 9, 11, 12, 13, 14, 16, 17, 22, 23, 24, 29. Moreover, for $n=19$ (resp. $n=28$), one of the two words is palindromic (resp. anti-palindromic), the other not. However, this phenomenon seem to disappear for larger $n$.
\item The value of $\minswentropy^{(2)}(n)/n$ increases with $n$ in general, but with the exceptions of $n=3, 4$, $n=6, 7$ and $n=12, 13$ (although the rounded numbers are the same). We believe that the exceptions are due to the effect of small size, and should not reproduce for larger $n$.
\item The number of runs for words of length $n$ achieving $\minswentropy^{(2)}(n)$ is weakly increasing with $n$, with the exception of $n=17, 28, 30, 35$. Moreover, for $n=28$ only one word among the two that has less runs than the word for $n=27$.
\item The maximal run length for words of length $n$ achieving $\minswentropy^{(2)}(n)$ is at most $3$, with the exception of $n=17, 28, 30, 31, 33, 37$. Moreover, for $n=28$, one of the two words has maximal run length $3$, and the other $4$.
\item There is only one word of length $n$ up to symmetries in \Cref{lem:switch} that achieves $\minswentropy^{(2)}(n)$, with the exception of $n=3, 19, 28$, where there are two such words.
\item Most words have only a unique most frequent subword satisfying the conditions in the caption of \Cref{tab:exhaustive-subwords}, with exceptions on $n = 5, 7, 8, 13, 14, 15, 17, 18$. Furthermore, for $n=19$, there is one word with only one most frequent subwords, but the other has three. We note that these exceptions happen at relatively small values of $n$.
\item There are roughly the same number of $0$'s and $1$'s in the words. Intuitively, when there is a significant difference in the number of symbols, subwords favoring the more frequent symbol will occur more often, raising the subword entropy. For instance, for words in which $0$ and $1$ occur with different frequency, \Cref{prop:maxocc-lowerbound} gives greater lower bounds.
\end{itemize}

We note that, for words of length $n$ achieving $\minswentropy^{(2)}(n)$, the ratio between the length of most frequent subwords and $n$ seems to stabilize around $0.4$, which is clearly different from the ratio $1/3$ for the lower bound given in \Cref{prop:maxocc-lowerbound}. It means that words achieving $\minswentropy^{(2)}(n)$ seem to have some kind of structure, in which there are few most frequent subwords. The bound by expectation thus should not be tight. Combined with the observation of \Cref{tab:exhaustive} that $\maxocc(w)$ for $w$ of length $n$ achieving $\minswentropy^{(2)}(n)$ has an exponential growth in $n$ with a growth constant close to but slightly larger than the value $2/3$ given by the lower bound in \Cref{prop:maxocc-lowerbound}. We thus have the following conjecture.

\begin{conj} \label{conj:limit-value-2}
  We have $L_2 > \logb(3/2)$.
\end{conj}

\Cref{conj:limit-value-2} is also supported by the fact that there are very few binary words of length $n$ with subword entropy close to the minimal $\minswentropy^{(2)}(n)$. Most words have clearly higher subword entropy, which implies that the lower bound in \Cref{prop:maxocc-lowerbound} valid for all words is improbable to be optimal. \Cref{conj:limit-value-2} implies that, for all large values of $n$, most frequent subwords of words achieving $\minswentropy^{(2)}(n)$ should also have its own structure, probably a high degree of self-correlation. The question remains on how to find such subwords.

However, we should note that we only have very limited data, as we were only able to perform exhaustive search for small values of $n$. A naïve method requires looking at $\Theta(4^n)$ word-subword pairs. Although some optimizations are possible, such as using \Cref{lem:switch} to reduce the number of words to examine, the time taken remains exponential, against which we cannot push too far. An evidence is that, although asymptotically $\minswentropy^{(2)}(n) / n$ should be bounded from below by $\logb(3/2) \approx 0.585$ by \Cref{coro:sw-bounds}, all the values of $\minswentropy^{(2)}(n) / n$ in \Cref{tab:exhaustive} are clearly smaller than this asymptotic bound, meaning that the values of $n$ tested here are not large enough. Nevertheless, we can still formulate reasonable conjectures based on these observations.

\begin{conj}
  For $k \geq 2$, let $w$ be a word of length $n \geq 1$ achieving the minimal subword entropy $\minswentropy^{(k)}(n)$. Then, except for a finite number of $n$, the longest run in $w$ has length $3$. Furthermore, the average run length converges when $n \to +\infty$.
\end{conj}

Intuitively, long runs in a word lead to higher subword entropy. Suppose that we have a run of length $k$ in a word $w$, for instance $w = w' \cdot 0^k \cdot w''$. Then by \Cref{lem:basics-3}, for any word $u', u''$, we have
\[
  \maxocc(w) \geq \occ(w, u' \cdot 0^{\lfloor k/2 \rfloor} \cdot u'') \geq \occ(w', u') \binom{k}{\lfloor k/2 \rfloor} \occ(w'', u''),
\]
However, $\binom{k}{\lfloor k/2 \rfloor}$ is the largest number of subword occurrences for words of length $k$. When $k$ is large, such a long run will thus lead to larger subword entropy, meaning that words achieving minimal subword entropy should not contain long runs.

\begin{conj}
  There is a finite number of $n$ such that there is a binary word $w$ of length $n$ achieving the minimal subword entropy with several most frequent subwords.
\end{conj}

We observe in \Cref{tab:exhaustive} that words achieving $\minswentropy^{(2)}(n)$ contains mostly runs of length $1$, $2$ and $3$. There are also longer runs, but they seem exceptional. It is this observation that leaded us to consider the three families of periodic words in \Cref{sec:binary} for better upper bound of $L_2$. Inspired by \Cref{thm:periodic-family}, we also looked experimentally at some other periodic words containing runs of length $1$, $2$ and $3$. Some of them seem to have the potential to give a slightly better upper bound of $L_2$. For instance, words of period $0001100111$ seem to have one most frequent subword that is periodic with period $0011$. The corresponding generating function is
\[
  f_{0001100111, 0011} = \frac{(1-x)^3 - x(9x^2 + 78x + 13)y}{(1-x)^4 - 9x(1-6x)^2y^2 - x(9x+16)(21x+4)y}.
\]
The largest coefficient in $[x^n]f_{0001100111, 0011}$ occurs asymptotically at the term $y^{\alpha n}$ with $\alpha \approx 1.004849\ldots$, whose expression involves the smallest real root of $z^4 - 517z^3 - 258z^2 - 77z + 1 = 0$. If we indeed have the correct structure of most frequent subwords of $(0001100111)^m$, then it gives an upper bound of $L_2$ that is roughly $0.63272\ldots$, improving the result from \Cref{thm:upper-bound}. However, we don't know how to prove the observed structure on most frequent subwords of such periodic words in this case and in general, which leads us to the following conjecture.

\begin{conj} \label{conj:periodic-subword-struct}
  For a given word $w$, there is a word $v$ such that, for all $m$ large enough, there is a most frequent subword of $w^m$ that takes the form $u \cdot v^r \cdot u'$, with $u$ and $u'$ of lengths bounded by $|v|$.
\end{conj}

If \Cref{conj:periodic-subword-struct} holds, then using arguments similar to those in \Cref{thm:upper-bound}, we can reduce the computation of $\maxocc(w^m)$ to that of maximizing $\occ(w^m, v^r)$ while losing only a multiplicative constant. We can then apply \Cref{thm:periodic-family} and tools in analytic combinatorics in several variables to obtain a better upper bound for $L_k$. Again, \Cref{conj:periodic-subword-struct} seems natural, intuitive and supported by experimental evidence, but we don't see how to settle it.

To solve these conjectures, we need a better understanding of subword occurrences in words. However, many intuitive ideas seem very difficult to prove. For instance, we have the following conjecture that is surprisingly difficult to tackle.

\begin{conj} \label{conj:increasing-entropy}
  For fixed $k \geq 2$, there is a value $N$ such that the function $\minswentropy^{(k)}(n)/n$ is increasing for $n \geq N$.
\end{conj}

We already know experimentally that \Cref{conj:increasing-entropy} is not universally true for all $n$ (see \Cref{tab:exhaustive}). This point needs to be addressed in possible proofs.

Another intuitive idea on most frequent subwords of a given word $w$ is that their length should be smaller than $|w|/2$. The reasoning is that longer subwords have more letters, thus more possible occurrences, but this effect only works up till length $|w|/2$. However, even such an intuitive idea, supported by \Cref{prop:maxocc-lowerbound}, seems difficult to prove.

\begin{conj} \label{conj:half-length}
  For a given word $w$ of length at least $2$, there is at least one most frequent subword of $w$ with length at most $\lceil |w| / 2 \rceil$.
\end{conj}

Should \Cref{conj:half-length} be correct, to compute the subword entropy of $w$, we only need to check subwords of length at most $\lceil |w| / 2 \rceil$, which would greatly speed up the current algorithm and allows us to reach higher values of $n$.

\bibliographystyle{alpha}
\bibliography{subword-fang}

\end{document}